\documentclass[letterpaper, 10 pt, conference]{ieeeconf}  

\IEEEoverridecommandlockouts                              

\usepackage{epsfig}
\usepackage{latexsym}
\usepackage{pdfsync}
\usepackage{amsmath,amssymb,amsfonts}

\def\tr{^{\rm T}}
\def\Real{{I\!\!R}}

\newcommand{\eq}{\begin{equation}\begin{array}{rllllllllllllllllllllllllllllllll}}
\newcommand{\ee}{\end{array}\end{equation}}
\newcommand{\bmt}{\left[ \begin{array}{ccccccccc}}
\newcommand{\emt}{\end{array}\right]}
\newcommand{\bea}{\begin{eqnarray}}
\newcommand{\eea}{\end{eqnarray}}
\newcommand{\bean}{\begin{eqnarray*}}
\newcommand{\eean}{\end{eqnarray*}}
\newcommand{\bc}{\begin{center}}
\newcommand{\ec}{\end{center}}

\title{\LARGE \bf
Linear Quadratic Regulation  for First Order Hyperbolic PDEs*
}

\author{Arthur J. Krener$^{1}$
\thanks{*This work was supported by AFOSR under grant  FA9550-20-1-0318}
\thanks{$^{1}$Arthur J. Krener is with the Department of Mathematics, University of California, Davis, CA 95616, USA
        {\tt\small ajkrener@ucdavis.edu}}%
}

\begin{document}
\maketitle

\begin{abstract}
We consider transport processes that are modeled by  first order hyperbolic partial differential equations.    Our goal is to find a full state feedback that makes a given reference profile locally asymptotically stable.  To accomplish this we employ Linear Quadratic Regulation (LQR) with finite dimensional patch or  point control actuation.    We derive the Riccati partial differential equation whose solution is the  kernel of the optimal cost. The optimal state feedback is also found. The derivation is accomplished by elementary techniques such as integration by parts and completing the square.  We apply this theory to two examples  that have appeared in the literature and that were solved by a modification of  LQR.   The first example deals with a model of a fixed-bed chemical reactor and the second example deals with traffic congestion on a stretch of freeway.
 \end{abstract}

\section{Introduction}
Aksikas et al. \cite{AFFW09} and Block and Stockar  \cite{BS24}, used a modification of infinite horizon Linear Quadratic Regulation (LQR)
to asymptotically stabilize to a reference profile  a system described by a first order hyperbolic partial differential equation.  Aksikas et al.  considered 
a model of a fixed-bed chemical reactor and Block and Stockar considered a model of traffic congestion on a length of freeway.   Both papers took a novel approach to solving the differential Riccati equation that arise in infinite dimensional LQR, equation (6.56) of \cite{CZ95}.   They assumed that the unknowns in 
the differential Riccati equation are not only the operator that is the kernel of the optimal cost but also the state and control weighting matrices that are in the criterion to be minimized.  With these extra degrees of freedom they were able to arrive at solutions to the differential Riccati equations.

The purpose of this paper is to offer a different approach where the only unknown in the differential Riccati equation is the kernel of the optimal cost.  And our differential Riccati equation is a partial rather than ordinary differential equation.
Both Aksikas et al. and Block and Stockar  assume  infinite dimensional distributed control actuation which is impossible to achieve.  We use only  finite dimensional patch or point control actuation which is more realistic.

The rest of the paper is organized as follows.  In the next section we introduce systems modeled by first order hyperbolic PDEs under finite dimensional patch or point control actuation.  In Section III we use infinte horizon LQR to derive the appropriate Riccati partial differential equation whose solution yields  a state feedback control law that locally, asymptoticaly stabilizes a reference profile of such a system. The next section applies these techniques to a model of a fixed-bed reactor considered by Christofides \cite{Ch01} and Aksikas et al. \cite{AFFW09}.  In Section V we apply these techniques to a model of traffic congestion similar to that considered by Block and Stokar \cite{BS24}.

\section{First Order Hyperbolic Systems }
Consider the linear, first order hyperbolic partial differential equation
\bea \label{dyn1}
\frac{\partial\zeta }{\partial t}(t,x)&=& D(x,\zeta)\frac{\partial \zeta }{\partial x}(t,x)+E(x,\zeta)z (t,x)
\eea
where $\zeta (t,x) \in \Real^n$, $t\in [0,T]$, $x\in [0,1]$,  $D(x,\zeta)$ is an $n\times n $ diagonal matrix valued function and $E(x,\zeta) $ is an $n\times n $ matrix.     The dynamics is subject to a given   initial condition $\zeta (0,x)$ and a given boundary condition  $\zeta (t,0)$. 
The control objective is to stabilize the dynamics to a given profile $\zeta^0(x)$.   We use the term "profile" for functions of $x$ alone and  the term "trajectory" for functions of  $t$ and $x$.
We  apply infinite horizon Linear Quadratic Regulation (LQR)  to the infinite dimensional system linearized about a reference profile.
  
We  consider two types of control actuation.  
Under patch actuation there are $m$ intervals $[a_j,b_j]$ with disjoint interiors. The controlled dynamics  is
\bea \nonumber
&&\frac{\partial \zeta }{\partial t}(t,x)= D(x,\zeta)\frac{\partial \zeta }{\partial x}(t,x)+E(x,\zeta)\zeta (t,x)\\
&&+G_j(x,\zeta)  \nu _j(t) \label{nldyn}
\eea
where $ G_j(x,\zeta)=G_j(\zeta)\chi_j(x)$,  $G_j(\zeta)$ is an  $n\times 1$ vector and $\chi_j(x)$ is the characteristic function of $ [a_j,b_j]$.
The $j^{th}$ control is denoted by $\nu_j(t)$.
Throughout this paper we invoke the convention of summing on repeated indices.  Typically an index $i$ runs from $1$ to $n$ while an index $j$ runs from $1$ to $m$. 

The other type  is
  point  actuation at locations $0\le a_1< a_2< \cdots < a_m\le 1 $
where   $G_j(x,\zeta)=G_j(\zeta)\delta(x-a^-_j)$.  Notice that the jumps take place on the left side of the $a_j$'s so
the resulting trajectories are continuous in $x$ from the right. 

We consider the linear approximating system around a reference state profile  $\zeta^0(x)$ generated by a reference control  trajectory $\nu^0(t)$.
We define displacement coordinates
\bean
z(t,x)= \zeta(t,x)-\zeta^0(x),&& u(t)=\nu(t) -\nu^0(t)
\eean
The linear approximating system is 
\bea \nonumber
&&\frac{\partial z }{\partial t}(t,x)= D^0(x)\frac{\partial z }{\partial x}(t,x)+E^0(x)z(t,x)\\
&&+G^0_j(x)  u _j(t) \label{ldyn}
\eea
where
$D^0(x)=D(x,\zeta^0(x))$, $E^0(x)=E(x,\zeta^0(x))$, $G_j^0(x)=G_j(x,\zeta^0(x))$.

\section{ Infinite Horizon LQR}
We consider infinite horizon LQR for the infinite dimensional linear system  (\ref{ldyn}).
  For background on  time invariant LQRs over an infinite horizon  for systems described by  PDEs we refer the reader to treatise \cite{CZ95}. 
Since  $D^0(x),\ E^0(x), \ G^0_j(x)$ are independent of $t$ we can employ infinite horizon LQR to find a stabilizing, time invariant feedback.   We choose  an $n\times n$ nonnegative definite matrix valued function  $Q(x_1,x_2)\ge 0$ that is symmetric in its arguments  $Q(x_1,x_2)= Q(x_2,x_1)$  and  an
$m\times m$ positive definite matrix $  R>0 $.  We   seek to minimize by choice of control $u(t)$ the criterion
\bea \label{criti}
&&{1\over 2}\int_0^\infty \iint_{[0,1]^2} z\tr(t,x_1)Q(x_1,x_2)z(t,x_2)\ dA\\
&&+u\tr(t)Ru(t)\ dt \nonumber
\eea
\normalsize
where $dA=dx_1dx_2$ and
 subject to the linear dynamics (\ref{ldyn}), a given initial condition $z(0,x)$ and a given boundary condition $z(t,0)$.

Let $P(x_1,x_2)$ be an $n\times n$ matrix valued $C^1$ function that is symmetric in its arguments $P(x_1,x_2)=P(x_2,x_1)$
and satisfies the boundary conditions $P(1,x_2)=P(x_1,1)=0$.   Assume that for every initial condition $z(0,x)$  and every boundary condition $z(t,0)$ there is a control trajectory $u(t)$ such that $z(t,x) \to 0$ as $t\to \infty$.  This is the stabilizability assumption that is required by LQR.  The required detectability 
condition is that if $z(t,x)$ is a trajectory sich that 
\bean
0&=&\iint_{[0,1]^2} z\tr(t,x_1)Q(x_1,x_2)z(t,x_2)\ dA
\eean
then $z(t,x)\to 0$  as $t\to\infty$. 

By the Fundamental Theorem of Calculus
\bean
&&0=  \iint_{[0,1]^2} z\tr(0,x_1)P(x_1,x_2)z(0,x_2)\ dA\\
&&+\int_0^\infty {d\over dt} \iint_{[0,1]^2} z\tr(t,x_1)P(x_1,x_2)z(t,x_2)\ dA \ dt
\eean

We bring the time differentiation inside the spatial integrals and integrate by parts to obtain
\small
\bean
&&0= \iint_{[0,1]^2} z\tr(0,x_1)P(x_1,x_2) z(0,x_2) \ dx_1dx_2\\
\\&&-\int_0^\infty\iint_{[0,1]^2}z\tr(t,x_1) (D^0)\tr(x_1)\frac{\partial P}{\partial x_1}(x_1,x_2)z(x_2)\ dA\ dt\\
&&-\int_0^\infty\iint_{[0,1]^2}z\tr(t,x_1) \frac{\partial (D^0)\tr}{\partial x_1}(x_1)P(x_1,x_2)  z(x_2) \ dA\ dt\\
&&-\int_0^\infty\int_0^1 z\tr(t,0) (D^0)\tr(0)P(0,x_2) z(t,x_2)\ dx_2\ dt\\
&& +\int_0^\infty\iint_{[0,1]^2}z\tr(t,x_1)(E^0)\tr(x_1) P(x_1,x_2)z(t,x_2) \ dA\ dt\\
&&-\int_0^\infty\iint_{[0,1]^2}z\tr(t,x_1)\frac{\partial P}{\partial x_2}(x_1,x_2)D^0(x_2)z(t,x_2) \ dA\ dt\\
&&-\int_0^\infty\iint_{[0,1]^2}z\tr(t,x_1)P(x_1,x_2)\frac{\partial D^0}{\partial x_2}(x_2)  z(t,x_2)\ dA\ dt\\
&&-\int_0^\infty\int_0^1z\tr(t,x_1)P(x_1,0)  D^0(0)z(t,0)\ dx_1\ dt\\
&& +\int_0^\infty\iint_{[0,1]^2}z\tr(t,x_1) P(x_1,x_2) E^0(x_2) z(t,x_2) \ dA\ dt\\
&&+  \int_0^\infty \iint_{[0,1]^2}  \left( G^0_j(x_1) u_j(t)\right)\tr P(x_1,x_2) z(t,x_2)\ dA\ dt\\
&&+  \int_0^\infty \iint_{[0,1]^2} z\tr(t,x_1)P(x_1,x_2) \left(G^0_j(x_2) u_j(t)\right) \ dA\ dt
\eean
\normalsize

We add the right side of this equation to the criterion (\ref{criti}) to get an equivalent criterion. Then we complete the square.  We seek a $K(x)$ such that the integrand of the time integral  of the equivalent criterion equals
\bean
&&\left(u(t) -\int_0^1 K(x_1)z(t,x_1)\ dx_1\right)\tr R\\
 &&\times \left(u(t) -\int_0^1 K(x_2)z(t,x_2)\ dx_2\right)\
\eean
The terms quadratic in $u(t)$ match so
equating terms bilinear in $u(t)$ and $z(t,x_2)$ we arrive at the formula 
\bean
K(x_2)&=& -R^{-1} \int_0^1(G^0)\tr(x_1)P(x_1,x_2)\ dx_1
\eean

Equating terms bilinear in $z(t,x_1)$ and $z(t,x_2)$ yields the  Riccati partial  differential equation 
\bea  \nonumber 
&&-(D^0)\tr(x_1)\frac{\partial P}{\partial x_1}(x_1,x_2)- \frac{\partial (D^0)\tr}{\partial x_1}(x_1)P(x_1,x_2)\\
&& \nonumber-\frac{\partial P}{\partial x_2}(x_1,x_2)D^0(x_2)-P(x_1,x_2) \frac{\partial D^0}{\partial x_2}(x_2)\\
&& \nonumber+\left(E^0(x_1)-\delta(x_1)D^0(x_1)\right)\tr P(x_1,x_2)\\
&&+P(x_1,x_2)\left(E^0(x_2)-D^0(x_2)\delta(x_2)\right)\nonumber \\
&& +Q(x_1,x_2) \nonumber  \\
&& \nonumber = \iint_{[0,1]^2}P(x_1,x_4)G^0(x_4) R^{-1} \\
&& \times (G^0)\tr(x_3)P(x_3,x_2)\ dx_3dx_4 \label{RPD} 
\eea
If   $P(x_1,x_2)$ and $ K(x_2)$ satisfy  these equations then the optimal cost starting at $z(0,x)$ is
\bean
 \iint_{[0,1]^2} z\tr(0,x_1)P(x_1,x_2) z(0,x_2) \ dx_1dx_2
\eean
and the optimal feedback is
\bean
u(t)&=&\int_0^1 K(x_2)z(t,x_2)\ dx_2
\eean

 \section{Example One:  Fixed-bed Reactor} 
 
 We consider a model of a fixed-bed reactor similar to those  studied by Christofides \cite{Ch01} and Aksikas et al. \cite{AFFW09}.   The  reactor   extends from $x=0$ to $x=1$ and the feed
species A is converted to product species B as it moves through the reactor.   The first state variable, ${\cal T}(t,x)$, 
is the temperature in the interior of the reactor and and the second state variable, ${\cal C}_A(t,x)$,  is the concentration of the reactant  A.  
Following \cite{AFFW09} we normalize these quantities
\bean
\theta_1(t,x)={{\cal T}(t,x)-{\cal T}_{in}\over {\cal T}_{in}},&& \theta_2(t,x)={{\cal C}_{A,in}-{\cal C}_A(x,t)\over {\cal C}_{A,in}}
\eean 

We assume that the reactor jacket is divided into  five sections  and we can control the temperature  independently on each section so this is an example of patch control.
The nonlinear dynamics is given by 
  \bea \label{nldyn}
&&  \frac{\partial \theta_1}{\partial t}(t,x)= -v_1 \frac{\partial \theta_1}{\partial x}(t,x)\\
&& \nonumber +k_1(1-\theta_2(t,x))\exp\left({ \nu  \theta_1(t,x)\over  \theta_1(t,x)+1}\right)\\
  &&
 + \beta(\nu_j(t)-\theta_1(t,x))  \nonumber \\
  && \frac{\partial \theta_2}{\partial t}(t,x)= -v_2\frac{\partial \theta_2}{\partial x}(t,x)  \nonumber \\
  && \nonumber +k_2(1-\theta_2(t,x))\exp\left({ \nu  \theta_1(t,x)\over  \theta_1(t,x)+1}\right)  \nonumber
 \eea
 where the control is ${\cal T}_j$, 
the temperature  of the $j^{th}$ section of the  jacket, $x\in [{j-1\over5}, {j\over 5}]$.  The controls are normalized to $\nu_j(t)= {{\cal T}_j-{\cal T}_{in}\over {\cal T}_{in}}$.  The constants $v_1,v_2,k_1,k_2.\nu ,\beta_j$ are given in \cite{AFFW09} between displays (19) and (20). Notice that 
 this is a quasi-linear PDE as the coefficients of the highest partial derivatives are constants.

 In contrast to our control on each of five patches, Christofides \cite{Ch01} assumes patch control with only one patch, the whole jacket of the reactor.    So his control is one dimensional.  Since $0<v_1<<v_2$ there are two time scales in the reactor.
 Christofides  exploits this by using singular perturbation techniques to  develop a  control law that   is robust to parameter uncertainties.  

  Aksikas et al.~\cite{AFFW09} use a modified version of infinite horizon LQR to develop a similar control law. They assume distributed control actuation, that is,  at each location the jacket temperature is independently controllable, so their control is infinite dimensional.  They solve a  "modified  Riccati ordinary differential equation" to find the linear feedback.  We call it "modified" because the unknowns in their equation are not only the kernel of the optimal feedback (in their terminology, the solution of the operator Riccati equation) but also the state  and control  weighting matrices in the running cost, (in standard notation, which they do not use, these are the $Q$ and $R$ matrices).  Aksikas et al. claim in Comment 3.1 (1) that their  "Riccati ordinary differential equation" is solvable based on a simple dimension count of the number of unknowns  vs the number of equations.    
  
   In Theorem 7 they assume that the kernel of the optimal cost, the state  and control  weighting matrices  are  diagonal  and they arrive at their equations (25).   These are two ordinary differential equations with quadratic nonlinearities and a quadratic side constraint.  They are somewhat vague about  what are the unkowns in these equations but they seem to imply that they are five functions,  the diagonal entries $\phi_1$ and $\phi_2$ of the kernel of the optimal cost and the diagonal enties $w_1$, $w_2$  of the state   weighting matrix.  The fifth unknown $p$ is redundant because it can  be easily absorbed into $w_1$ and $w_2$.  
  They do not prove that these equations have a nonnegative solution but in Figure 2 they present a numerically computed $\phi_1$ that is nonnegative.

We   use infinite horizon, infinite dimensional LQR to find a linear feedback
  that robustly stabilizes a fixed bed reactor to a given profile  $z^0(x)$.   Because the refence profile does not depend  on time, the linearization around it is time invariant. We consider that the partial differential equation (\ref{RPD}) to be the   Riccati equation for this problem.  
 We choose the a similar reference profile $\theta^0(x)$  as \cite{AFFW09} by setting  the same initial conditions as  \cite{AFFW09}, ${\cal T}_{in}=320$ K, ${\cal C}_{A,in}= 4$.  But we use patch control with five patches instead of the distributed control of \cite{AFFW09}.  The temperture of the $j^{th}$ patch is denoted by ${\cal T}_j$ for $j=1,\ldots,5$.    We denote the normalized reference controls by $\nu^0_j = {{\cal T}_{j}^0\over 320}={5\over 4}$.
   
 We use the Method of Lines to to approximately solve the nonlinear PDE (\ref{nldyn}) for the reference profile.  We choose an $N$ and subdivide $[0,1]$ into $N$ subintervals with endpoints $\xi_k={k\over N},\ k=0,\ldots,N$. We approximate $\theta^0_{i}(\xi_k)$ by $\theta^0_{i,k}$.  We set $\theta^0_{i,0}=0$ and use backward spatial differences to approximate the  spatial derivatives for $k=1,\ldots,N$.
 This leads to the $2N$ dimensional nonlinear ODE
\bean
\dot{\theta}^0_{1,k}(t)&=&-v_1N\left(\theta^0_{1,k}(t)-\theta^0_{1,k-1}(t)\right)\\
&&+k_1\left(1-\theta^0_{2,k}(t)\right)
\exp\left({\nu  \theta^0_{1,k}(t)\over \theta^0_{1,k}(t)+1}\right)\\
&&+\beta\left(400-\theta^0_{1,k}(t)\right)\\
\dot{\theta}^0_{2,k}(t)&=&-v_2N\left(\theta^0_{2,k}(t)-\theta^0_{2,k-1}(t)\right)\\
&&+k_2\left(1-\theta^0_{2,k}(t)\right)
\exp\left({\nu  \theta^0_{1,k}(t)\over \theta^0_{1,k}(t)+1}\right)
\eean
for $k=1,\ldots,N$.

 The resulting reference state profile in the original ${\cal T}^0(t,x),\  {\cal C}^0(t,x)$ coordinates is seen in Figures 
1 and 2 with $N=100$. These  profiles are different from  Figure 1 of \cite{AFFW09} because we have chosen different constants.
 
 \begin{figure}[h]
\centering
\begin{minipage}{.45\linewidth}
  \includegraphics[width=\linewidth]{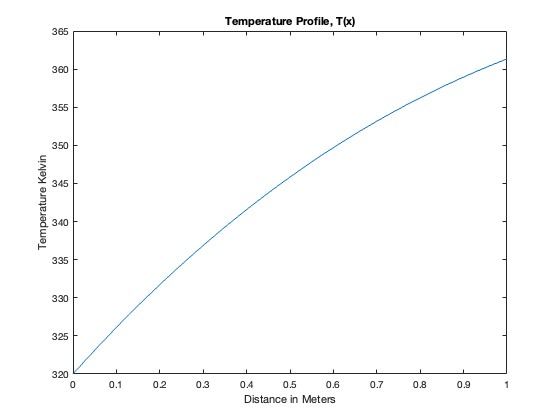}
  \label{TP}
  \caption{Temperature  Profile }
\end{minipage}
\hspace{.05\linewidth}
\begin{minipage}{.45\linewidth}
  \includegraphics[width=\linewidth]{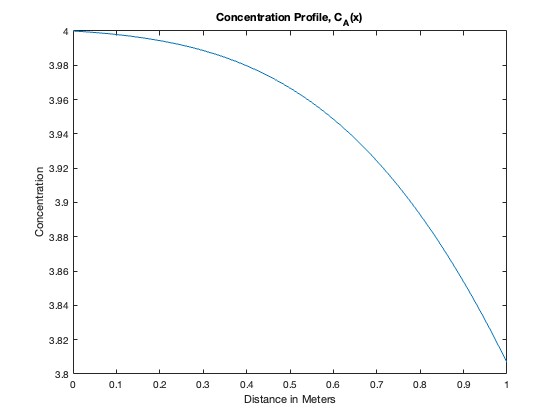}
  \label{CP}
  \caption{Concentration Profile }
  \end{minipage}
\end{figure}

 We define  perturbation variables
 \bean
 z_i(t,x)=\theta_i(t,x)-\theta_i^0(x),&&   u_j(t)= \nu_j(t)-\nu_j^0
\eean 
 then 
 the perturbation dynamics through linear terms is given by  
 \bea \label{perdyn}
 && \frac{\partial z}{\partial t}(t,x)=D^0 \frac{\partial z}{\partial x}(t,x) + E^0(x)z(t,x)\\
  && \nonumber +G^0_j(x)u_j(t)
   \eea
 where 
 \bea \label{notation}
D^0&=& \bmt -v_1&0\\0&-v_2\emt\\  \nonumber
E^0(x)&=&\exp(\mu ^0(x))\bmt k_1\mu {1-\theta^0_{2}(x)\over (1+\theta^0_{1}(x))^2} &-k_1\\
\\ k_2\mu {1-\theta^0_{2}(x)\over (1+\theta^0_{1}(x))^2} &-k_2\emt \\ \nonumber
G^0_j(x)&=& \beta \chi_j(x)
\eea
for $j=1,\ldots,5$
where $\chi_j(x)$ is the characteristic function of $ [0.2(j-1),0.2j] $
 and $\mu ^0(x)={\mu  \theta^0_{1}(x)\over 1+\theta^0_{1}(x)}$.
Aside from a change of notation this is the same as
 (22) of \cite{AFFW09}.
  
  We discretize the perturbation dynamics (\ref{perdyn}) by  letting $\zeta_k(t) \approx z(t, \xi_k)$,
 and  $\frac{\partial z}{\partial x}(t,\xi_k)\approx N(\zeta_{k+1}(t) -\zeta_k(t) $ for $k=0,\ldots, N$.
 We choose $Q(x_1,x_2) =I^{2\times2}\delta(x_1-x_2)$ which we discretize by $Q_{k_1,k_2}=\delta_{k_1,k_2}$ and $R =I^{5\times5}$.  Then we use MATLAB's 
 {\tt lqr.m} to approximately solve the LQR.  The $1-1$ and $2-2$   blocks of the kernel $P$ of the optimal cost are shown in Figures 3 and 4.  The least stable
 closed loop eigenvalues are $-7968\pm 3204 i$  so the system is quite stable.   
  Figures 3 and 4 show the $1,1$ and $2,2$ components of $P$.
  \begin{figure}[h]
\centering
\begin{minipage}{.45\linewidth}
  \includegraphics[width=\linewidth]{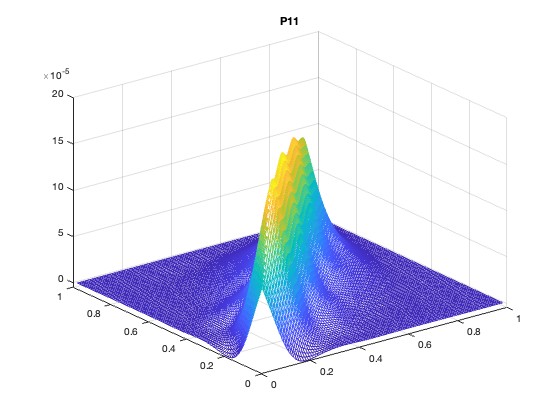}
  \label{P11}
  \caption{$P_{1,1}$ }
\end{minipage}
\hspace{.05\linewidth}
\begin{minipage}{.45\linewidth}
  \includegraphics[width=\linewidth]{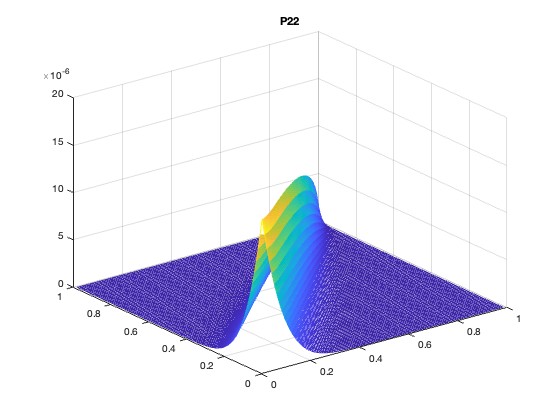}
  \label{P22}
  \caption{$P_{2,2}$ }
  \end{minipage}
\end{figure}

\section{Example Two: Traffic Control}

 Block and Stockar  \cite{BS24}
used a similarly modified version of   infinite horizon LQR to mitigate traffic congestion on a stretch of freeway  via variable speed limits.
Their control is distributed, it is the rate of change of the  speed limit which they assume can be set independently  at each point of the freeway.  So their control is infinite dimensional.

We use  infinite horizon LQR to mitigate traffic congestion using
metered on-ramps so our control is finite dimensional, the  dimension of the control is the number of interchanges in the stretch of freeway under consideration.  

Consider a stretch of freeway of length $10$ kilometers  starting at $x=0$, ending at $x=10$ with interchanges  at $x=a_j=2j,\  j=1,\ldots,4$.  
At each interchange $a_1, \ldots a_4$ cars can enter or leave the freeway.  The four dimensional control is pointwise, the net rate of cars entering or leaving the freeway at each interchange. The goal is to find a strategy to meter
the rate of cars entering the freeway so as to avoid congestion.  

Like Block and Stockar we use the Lighthill-Whitam-Richards model (LWR) to describe the traffic flow.  
There are more sophisticated models such as the Aw-Rascle-Zhang model that is used in \cite{YK19} to control 
the traffic but the LWR model is simpler and widely used.
The LWR model is a first order hyperbolic PDE,
\bea
\label{pdet}
\frac{\partial \rho}{\partial t}(t,x) &=& -\frac{\partial q}{\partial x}(t,x)
\eea
where $\rho(t,x)$ is  density of traffic (cars per kilometer) and $q(t,x)$
is the equilibrium relationship between traffic density and traffic flux.
We use Greenshield's equilibrium relationship between typical traffic velocity $v$ (kilometers per minute) and typical  traffic density $\rho$ (cars per kilometer)
\bea \label{GR}
v=v_M\left(1-{\rho\over \rho_M}\right)
\eea
where $v_M$ is the maximum velocity and $ \rho_M$ is the maximum density.
This leads to a quadratic relationship between traffic flux $q$ and  traffic density
$\rho$,
\bean
q(\rho)&=&v_M \left(1- {\rho\over \rho_M}\right)\rho 
\eean

The traffic flux function $q(\rho)$ is an inverted parabola with zeros at $\rho=0$ and $ \rho=\rho_M$ and
a maximum flux at the critical density $ \rho_C={\rho_M\over 2}$.  The traffic flux is increasing for $0<\rho<\rho_C$
and deceasing  for $\rho_C<\rho<\rho_M$.  The critical density separates the free flow region  from 
the congested region. The critical velocity $v_C={v_M\over 2}$ is the typical velocity at the critical density according to Greenshield's relation (\ref{GR}).

Plugging Greenshield's relation (\ref{GR}) into (\ref{pdet}) we obtain the nonlinear PDE
\bea
\label{pdet1}
\frac{\partial \rho}{\partial t}(t,x) &=&-v_M \left( 1-{2\rho(t,x) \over \rho_M}\right) \frac{\partial \rho}{\partial x}(t,x)
\eea
Notice that this is not a quasi linear PDE because the coefficient of the spatial partial derivative of the dependent variable $\rho(t,x)$ contains the dependent variable.

The initial condition is a given $\rho(0,x)$ and the boundary condition, the number of cars  entering the freeway at its beginning $a_0=0$, is a given $\rho(t,0)$. 
  The controls are the net rate
 of cars  entering the freeway at $a_j=2j,\ j=1,\ldots,4$. 
\bea \label{in}
\frac{\partial  \rho}{\partial t}(t,a^+_j)-\frac{\partial  \rho}{\partial t}(t,a^-_j)&=& G_j \nu_j(t)
\eea
where units of  $G_j$ are inverse kilometers and the units of $\nu_j(t)$ are net rate of cars  entering the freeway.  
The goal is to stabilize the traffic to a the reference density profile $\rho^0(x)< \rho_C$ that satisfies (\ref{pdet1}).
 This is an example of point actuation. 
 
 We choose a
 reference profile $\rho^0(x)$ that is not time varying,
 $\frac{\partial \rho^0}{\partial t}(x)=0$ except at the interchanges $x=a_j=2j,\ j=1,\ldots,4$  where $\frac{\partial \rho^0}{\partial t}(x)=0$ can jump. Then (\ref{pdet1})
implies $\frac{\partial \rho^0}{\partial x}(x)=0$  except at the interchanges.
The reference profile  is defined by $\rho^0(0)=0.9\rho_C$ with jumps  at each of the interchanges by  $ \nu^0_j=0.02 \rho_C, \ j=1,\ldots,4$.  In other words
\bean
\rho^0(x)&=&(0.9+0.02j)\rho_C
\eean
for $x\in (a_j, a_{j+1})$, $j=0,\ldots,4$
where we have taken $G_j=1$.  The reference profile  is less  that the critical density 
peaking at $0.98 \rho_C$ in the last two kilometers.

We linearize the dynamics around this profile. We denote the actual trajectory by $\rho(t,x)$ and the actual controls by
 $\nu_j(t),\ j=1,\ldots,4$.  The actual number of cars entering at the beginning of the stretch of freeway is $\rho(t,0)$.  We define perturbation variables
\bean
z(t,x)= \rho(t,x) - \rho^0(x),&& u_j(t) =\nu_j(t)-\nu_j^0
\eean

The nonlinear model in perturbation variables is
\bea \label{dynnon}
\frac{\partial z}{\partial t}(t,x) =-v_M \left( 1-2{\rho^0(x)+z(t,x) \over \rho_M}\right) \frac{\partial z}{\partial x}(t,x)
\eea
subject to the boundary conditions
\bea  \label{zbc0}
z(t,0) &=& \rho(t,0)-0.9\rho_C\\
\frac{\partial z}{\partial t}(t,a^+_j)&=&\frac{\partial  z}{\partial t}(t,a^-_j)+u_j(t)  \label{zbcj}
\eea
 and the initial condition 
\bea \label{zic}
z(0,x)&=&\rho(0,x)-\rho^0(x)
\eea

We assume that $z(t,x)$ is small  enough so that 
we can neglect the term quadratic in $z(t,x)$
 and so  we obtain the linear model
  \bea \label{dynlin}
\frac{\partial z}{\partial t}(t,x) &=&-v_M \left( 1-{2 \rho^0(x) \over \rho_M}\right) \frac{\partial z}{\partial x}(t,x)
\eea
subject to the boundary conditions (\ref{zbc0}), (\ref{zbcj}) and the initial condition (\ref{zic}).

We discretize this model by choosing a integer $N>0$ divisible by $5$ and subdivide $[0,10]$ into $N$ subintervals with endpoints $\xi_k={10 k\over N}$ for 
$k=0,\ldots,N$.  We approximate $z(t,\xi^+_k)$ by $z_k(t)$ and $ \frac{\partial z}{\partial x}(t,xi_k)$ by ${N\over 10}(z_k(t) -z_{k-1}(t))$.  We use backward spatial differences because of the boundary condition on $z(t,0)=z_0(t)$.   The Method of Lines  leads to finite dimensional system
\bean
\frac{dz_k}{dt}(t)&=& -v_M \left( 1-{2 \rho^0(\xi_k) \over \rho_M}\right) {N\over 10}(z_k(t) -z_{k-1}(t))\\
&&+u_{j}(t)
\eean
if $k={N\over 5}j$ for $j=1,\ldots,4$ and
\bean
\frac{dz_k}{dt}(t)&=& -v_M \left( 1-{2 \rho^0(\xi_k) \over \rho_M}\right) {N\over 10}(z_k(t) -z_{k-1}(t))
\eean
othewise.

This is the time invariant linear system of the form
\bean
\dot{z}&=& Fz+Gu
\eean 
with state dimension $N$, control dimension $m=4$,
\bean
F_{k,k}&=& - {N\over 10}v_M \left( 1-{2 \rho^0(\xi_k) \over \rho_M}\right)\\
F_{k,k-1}&=&  {N\over 10}v_M \left( 1-{2 \rho^0(\xi_k) \over \rho_M} \right)\\
G_{k,j}&=& 1\quad\ \mbox{ if } k={N\over 5}j
\eean
and zero otherwise.

If we take $Q(x_1,x_2)=\delta(x_1-x_2)$, $R(t)=I^{4\times 4}$ and ${\cal P}(x_1,x_2)=\delta(x_1-x_2)$ then the criterion is
\bean
\int_0^\infty \int_0^{10} \left|z(t,x)\right|^2 \ dx + \left|u(t)\right|^2\ dt 
\eean
which we discretize to
\bean
\int_0^\infty {10\over N} \sum_{k=1}^N z^2_k(t) +\sum_{j=1}^4 u_j^2 \ dt
\eean

In other words the discrete versions of these are
\bean
Q_{k_1,k_2}= {10\over N}\delta_{k_1,k_2},&& R_{j_1,j_2}=\delta_{j_1,j_2}
\eean

We choose parameters similar to those in \cite{BS24},
\bean
\begin{array}{|c|c|c|c|c|}
\hline
 \mbox{Parameter}&\mbox{Symbol}& \mbox{Value}& \mbox{Unit}\\
 \hline
 \mbox{Maximum density}& \rho_M& 160& \mbox{cars/km}\\
 \mbox{Critical  density}& \rho_C& 80& \mbox{cars/km}\\
  \mbox{Maximum speed}&v_M& 2&\mbox{km/min}\\
  \mbox{Critical speed}&v_C&1&\mbox{km/min}\\
  \mbox{Road length}&L& 10&\mbox{km}\\
  \hline
  \end{array}
\eean
and $N=100$.

\begin{figure}[h] 
\centering
\includegraphics[width=3in]{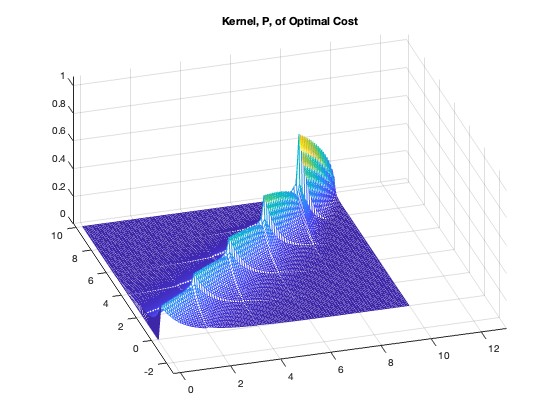} 
\caption{Kernel, $P$, of Optimal Cost}
\label{Pt}
\end{figure}

We use MATLAB's {\tt lqr.m} to solve the corrseponding Riccati equation.
The kernel $P$ of the optimal cost is shown in Figure \ref{Pt}.
The least stable closed loop eigenvalue is $-0.2567$.

\section{Conclusion}
We have  used Linear Quadratic Regulation to locally stabilize to a reference profile a system described
by a first order hyperbolic linear partial equation under two types of finite dimensional control actuation.
Patch actuation is when there are a finite number of intervals in the spatial domain and on each interval a control acts 
uniformly.  Point actuation is when there are a finite number of points in the spatial domain and a control acts 
at each point.  
   
   We applied our techniques to two hyperbolic problems that have appeared in the literature and were solved by a modification of LQR.
 The first problem is the stabilization of a fixed-bed chemical reactor to a given reference profile using control actuation on five patches.
The second problem is the stabilization of a stretch of freeway to a given reference profile using  control actuation at four points.

\end{document}